\documentclass[12pt]{article}
\usepackage[all]{xy}
\usepackage{amsmath}
\usepackage{amsfonts}
\usepackage{amssymb}
\usepackage{amscd}
\usepackage{amsthm}
\usepackage{latexsym}
\usepackage{amsbsy}

\newtheorem{prop}{Proposition}[section]
\newtheorem{theorem}{Theorem}[section]
\newtheorem{corollary}{Corollary}[section]

\begin{document}
\title
{ Cyclic Cohomology of Crossed Coproduct Coalgebras}
\author{R. Akbarpour,  M. Khalkhali \\Department of Mathematics 
, University of Western Ontario\\
\texttt{akbarpur@uwo.ca \; masoud@uwo.ca}}
\maketitle

\begin{abstract}
We extend our work in~\cite{rm01} to the case of Hopf comodule coalgebras.
We introduce the cocylindrical module $C \natural^{} \mathcal{H}$, where $\mathcal{H}$ is a Hopf algebra with
bijective antipode and $C$ is a Hopf comodule coalgebra over $\mathcal{H}$. We show that there exists an isomorphism between
the cocyclic module of the crossed coproduct coalgebra 
$C \! > \! \blacktriangleleft \! \mathcal{H} $ and $\Delta(C \natural^{} \mathcal{H}) $, the cocyclic 
module related to the diagonal of $C \natural^{} \mathcal{H}$. We approximate $HC^{\bullet}(C \!> \! \blacktriangleleft \! \mathcal{H}) $ by a spectral
sequence and we give an interpretation for $ \mathsf{E}^0 , \mathsf{E}^1$ and $\mathsf{E}^2 $ terms of this 
spectral sequence.  
\end{abstract}


\section{Introduction}
Getzler and Jones in~\cite{gj93} introduced a method to compute the cyclic homology of a crossed product algebras $ A \rtimes G $,
where $G$ is a group that acts on an algebra $A$ by automorphisms. This method is based on constructing a cylindrical module, $A \natural G$, and
showing that $ \Delta( A \natural G) \cong \mathsf{C}_{\bullet}
(A \rtimes G)$, where $\Delta$ is the diagonal and $\mathsf{C}_{\bullet}$ the cyclic module functor. Then by using 
the Eilenberg-Zilber theorem for cylindrical modules, they obtained a quasi-isomorphism of mixed 
complexes $\Delta(A \natural G) \cong Tot_{\bullet}{(A \natural G)}$, and a 
spectral sequence converging to $HC_{\bullet}(A \rtimes G)$. This spectral sequence was first obtained, by a different method, 
by Feigin and Tsygan~\cite{fs86}.
We used the the same method in~\cite{rm01} to generalize their work to Hopf module algebras with the action of a Hopf algebra
$\mathcal{H}$ on an algebra $A$, where the antipode of $\mathcal{H}$ is assumed to be bijective.\\
In this paper we continue our work in~\cite{rm01} to obtain similar results for Hopf comodule coalgebras.

We construct a cocylindrical module $C \natural^{} \mathcal{H}$ and we show that there exists an isomorphism
between cocyclic structures of the diagonal of $C \natural^{} \mathcal{H},$ denoted by $\Delta(C \natural^{} \mathcal{H})$,  
and the crossed coproduct coalgebra $ C \!>\! \blacktriangleleft \! \mathcal{H}$. We apply the Eilenberg-Zilber theorem for cocylindrical modules to get a 
quasi-isomorphism of mixed complexes,  $\Delta(
C \natural^{} \mathcal{H}) \cong Tot^{\bullet}(C \natural^{} \mathcal{H})$, and we approximate   
$HC^{\bullet}( C \!> \! \blacktriangleleft \! \mathcal{H})$ by a spectral sequence. We give an interpretation for the first
three terms of this spectral sequence. In our computations we find the coinvariant cocyclic module $C^{\bullet}_{\mathcal{H}^{}}(C)$
that is dual to $C_{\bullet}^{\mathcal{H}}(A)$ that we found in~\cite{rm01}.\\


 \section{Preliminaries on cocylindrical modules and Hopf algebras   
}

Let $[j],  j \in \mathbb{Z}^{+} $ be the set of all integers with the standard ordering and 
$t_j$ an automorphism of $[j]$ defined by $ t_j(i)=i+1,i \in \mathbb{Z} $. The paracyclic
category $\Lambda_{\infty}$  is a category whose
objects are $[j], j \ge 0$ and whose morphisms, $Hom_{\Lambda_{\infty}}([i],[j])$, are the
sets of all nondecreasing maps $f:[i] \rightarrow [j]$ with $t_j^{j+1}f=f \; t_i^{i+1}$~\cite{fs86}.

One can identify the simplicial category $\Delta$ as a sub category of $\Lambda_{\infty}$ so that the objects 
of $\Delta$ are the same as $\Lambda_{\infty}$ and morphisms, $Hom_{\Delta}([i],[j])$, are those morphisms of 
$Hom_{\Lambda_{\infty}}([i],[j])$ that map $\{0,1, \dots ,i \}$ into $\{ 0,1, \dots , j \}$.\\
Every element $g \in Hom_{\Lambda_{\infty}}([i],[j]) \; \; , 0 \le r \le \infty ,$ has a unique
decomposition $g=f\; t_i^k $ where $ f \in Hom_{\Delta}([i],[j])$.

If $\mathcal{A}$ is any category, a paracocyclic object in $\mathcal{A}$ 
is a covariant functor from $\Lambda_{\infty}$ to $\mathcal{A}$. This definition is equivalent to giving 
a sequence of objects $A_0 , A_1 , \dots$ together with coface operators 
$\partial^i : A_n \rightarrow A_{n+1}, (i=0,1,\dots ,n+1),$ 
codegeneracy operators $\sigma^i : A_n \rightarrow A_{n-1}, (i=0,1,\dots ,n-1),\;$and cyclic operators 
$\tau_n:A_n \rightarrow A_n$ where these operators
satisfy the cosimplicial
and $\Lambda_{\infty}$-identities:      
\begin{eqnarray*} \label{eq:sim}
\partial^j \partial^i&= \partial^{i} \partial^{j-1}\hspace{15pt} i<j,\\ \notag
\sigma^j \sigma^i &= \sigma^{i} \sigma^{j+1} \hspace{15pt} i\le j,  \notag
\end{eqnarray*}
\begin{equation}
\begin{align}
\sigma^j \partial^i &=
\begin{cases}
 \partial^{i} \sigma^{j-1}&\text{$i<j$}\\ \notag
\mbox{identity}      &\text{$i=j$ or $i=j+1$}\\ \notag
\partial^{i-1} \sigma^j  &\text{$i>j+1$}.
\end{cases}
\end{align}
\end{equation}
\begin{equation} \label{eq:pa}
\begin{align}
\tau_{n+1} \partial^i = \partial^{i-1} \tau_{n},\hspace{15pt} 1\le i \le n\; , \quad
\tau_{n+1} \partial^0  = \partial^{n}, \notag \\
\tau_{n-1} \sigma^i  =  \sigma^{i-1} \tau_{n},  \hspace{15pt} 1\le i\le n\; , \quad
\tau_{n-1} \sigma^0  = \sigma^{n-1} \tau_{n}^2.  \\
\end{align}
\end{equation}
If in addition we have $\tau_n^{n+1}=id$, then we have a cocyclic object in the sense of Connes~\cite{aC994}.
By a bi-paracocyclic object in a category $\mathcal{A}$, we mean a paracocyclic object in
the category of paracocyclic objects in $\mathcal{A}$. So, giving a bi-paracocyclic object
in  $\mathcal{A}$ is equivalent to giving a double sequence $A(p,q)$ of objects of $\mathcal{A}$ and operators
$\partial_{p,q} ,\sigma_{p,q},\tau_{p,q}$ and
$ \bar{\partial}_{p,q} ,\bar{\sigma}_{p,q},\bar{\tau}_{p,q}$ such that, for all $p \ge 0$,
\begin{eqnarray*}
B_p(q)=\{ A(p,q),\sigma^i_{p,q} ,\partial^i_{p,q} ,\tau_{p,q} \},
\end{eqnarray*}
and for all $q \ge 0$,
\begin{eqnarray*}
\bar{B}_q(p)=\{ A(p,q),\bar{\sigma}^i_{p,q} ,\bar{\partial}^i_{p,q} ,\bar{\tau}_{p,q} \},
\end{eqnarray*}
are paracocyclic objects in $\mathcal{A}$ and every horizontal operator commutes with every vertical operator.

We say that a bi-paracocyclic object is cocylindrical~\cite{gj93} if for all $p,q \ge 0,$
\begin{equation} \label{eq:cy}
\bar{\tau}_{p,q}^{p+1} \;
\tau_{p,q}^{q+1}=id_{p,q}.
\end{equation}

If $A$ is a bi-paracocyclic object in a category $\mathcal{A}$, the paracocyclic
object related to the diagonal of $A$ will be denoted by $\Delta A$. So, the paracocyclic operators on $\Delta A(n)=A(n,n)$ are $\bar{\partial}^i_{n,n+1}
\; \partial^i_{n,n} , \bar{\sigma}^i_{n,n-1} 
\sigma^i_{n,n},\bar{\tau}_{n,n}  \tau_{n,n}$.
When $A$ is a cocylindrical object since the cyclic operator of $\Delta A $ is $\bar{\tau}_{n,n} \tau_{n,n}$ and  $\bar{\tau},\tau$
commute, then, from
$\bar{\tau}^{n+1}_{n,n} \tau^{n+1}_{n,n} =id_{n,n}$, we conclude that $(\bar{\tau}_{n,n} \tau_{n,n})^{n+1}=id$. So that $\Delta A$ is a cocyclic object.

Let $k$ be a commutative unital ring. By a paracocyclic (resp. cocylindrical or cocyclic) module over $k$, we mean a paracocyclic (resp. cocylindrical or cocyclic) object in the 
category of $k$-modules.\\

An important example of a cocyclic module in this paper is the cocyclic module of a coalgebra $C$, denoted by $\mathsf{C}^{\bullet}(C)$.
It is defined by $\mathsf{C}^n(C)=C^{\otimes(n+1)}, n \ge 0,$ with coface, codegeneracy and cyclic oparators 
\begin{eqnarray*}
& &\partial^i(a_0,a_1,\dots,a_n)=(a_0,\dots,a_i^{(0)},a_i^{(1)},\dots,a_n), \;\; 1 \le i \le n,\\
& &\partial^{n+1}(a_0,a_1,\dots,a_n)=(a_0^{(1)},\dots,\dots,a_n,a_0^{(0)}), \\
& &\sigma^i(a_0,a_1,\dots,a_n)=\epsilon(a_{i})(a_0,\dots,a_{i-1},a_{i+1},\dots,a_n), \;\; 1 \le i \le n,\\
& &\tau(a_0,a_1,\dots,a_n)=(a_1,a_2,\dots,a_0). 
\end{eqnarray*}

A paracochain complex, by definition, is a graded $k$-module $\mathsf{V}^{\bullet} = ( V^i)_{i \in \mathbb{N}}$\; equipped
with operators $b: V^i \rightarrow V^{i+1}$ and $ B: V^i \rightarrow V^{i-1} $ such that $ b^2 = B^2 =0 ,$ and the operator
$T = 1-(bB+Bb)$ is invertible. In the case that $T=1$, the paracochain complex is called a mixed (cochain) complex.

Corresponding to any paracocyclic module $A$, we can define the paracochain complex $\mathsf{C}^{\bullet}(A)$ with the underlying graded 
module   $\mathsf{C}^n(A)= A(n)$ and the operators $ b = \sum_{i=0}^n (-1)^i \partial^i$ and 
$ B =  N \sigma (1 - (-1)^{n+1} \tau ) $. Here, $\sigma$ is the extra degeneracy satisfying $\tau \sigma^0 = \sigma \tau $, and
 $ N= \sum_{i=0}^n(-1)^{in} \tau^i $ is the norm operator. For any bi-paracocyclic module $A$, $Tot(\mathsf{C}(A))$
is a paracochain complex with $Tot^n(\mathsf{C}(A))= \sum_{p+q=n} A(p,q)$ and with the operators $ Tot(b) = b + \bar{b} $ and
$Tot(B)=  B + T \bar{B}$, where $ T = 1-(bB+Bb)$. It is a mixed complex if $A$ is cocylindrical~\cite{gj93}.

If we define the normalized cochain functor $\mathsf{N}$ from  paracocyclic modules to paracochain complexes with the underlying
graded module $\mathsf{N}^n(A) = {\bigcap}_{i=0}^{n-1} ker(\sigma^i) $ and the operators $b ,B$ induced from $\mathsf{C}^{\bullet}(A)$, then we have the 
following well-known results(see~\cite{gj93} for a dual version):\\
1. The inclusion $ ( \mathsf{N}^{\bullet}(A) , b) \rightarrow (\mathsf{C}^{\bullet}(A),b)$ is a quasi-isomorphism of complexes.\\
2. The cyclic Eilenberg-Zilber theorem holds for cocylindrical modules, i.e.,  \\
for any cocylindrical module $A$, there is a natural quasi-isomorphism $\mathbf{f}_0 + \mathbf{u f_1} : \mathsf{N}^{\bullet}(\Delta(A))
 \rightarrow Tot^{\bullet}(\mathsf{N}(A))$ of mixed complexes, where $\mathbf{f}_0$ is the shuffle map.

When $A$ is a cocyclic module, then we can construct the related mixed complex $(\mathsf{C}^{\;\bullet}(A)\mathit{[}\mathbf{u}\mathit{]} , b + \mathbf{u} B)$.
The cohomology of the complex $(\mathsf{C}^{\;\bullet}(A)\mathit{[}\mathbf{u}\mathit{]}  \otimes_{k[\mathbf{u}]} \mathsf{W} , b + \mathbf{u} B)$ 
defines the cyclic cohomology of $A$ with coefficients in $\mathsf{W}$, where $\mathsf{W}$ is a graded 
$k[\mathbf{u}]$-module with finite homological dimention.
We denote $\mathsf{C}_{\;\bullet}(A) \otimes_{k[\mathbf{u}]} \mathsf{W} $ by  $ \mathsf{C}(A) \boxtimes
\mathsf{W} $. We know that if
 $\mathsf{W}=  k[\mathbf{u , u^{-1}}]$, $k[\mathbf{u , u^{-1}}]/\mathbf{u}k[\mathbf{u}],$\\
$k[\mathbf{u}]/\mathbf{u} k[\mathbf{u}]$ we get, respectively,  $HP^{\bullet}(A)$  the periodic cyclic cohomology, $HC^{\bullet}(A)$ the cyclic cohomology and $HH^{\bullet} (A)$ the Hochschild cohomology of $A$~\cite{gj93}.

In this paper the word coalgebra means a coassociative, not necessarily cocommutative, counital coalgebra over a fixed commutative ring $k$.
Similarly, our Hopf algebras are over $k$ and are not assumed to be commutative or cocommutative. 
The undecorated tensor product $\otimes$ means tensor product over $k$. If $\mathcal{H}$ is a Hopf algebra, we denote its coproduct
by $\Delta : \mathcal{H} \rightarrow \mathcal{H} \otimes \mathcal{H}$, its counit by $\epsilon : \mathcal{H} \rightarrow k$,
its unit by $\eta:k\rightarrow \mathcal{H}$ and its antipode by $S:\mathcal{H} \rightarrow \mathcal{H}$. We will use Sweedler's
notation $\Delta(h)=h^{(0)} \otimes h^{(1)}, \Delta^2(h)=(\Delta \otimes id)\circ \Delta (h)=(id \otimes \Delta )\circ \Delta(h)=h^{(0)} \otimes h^{(1)}\otimes h^{(2)}$, 
etc., where  summation is understood.\\
If $\mathcal{H}$ is a Hopf algebra, the word $\mathcal{H}$-comodule means a comodule over the underlying coalgebra of $\mathcal{H}$.
We use Sweedler's notation for comodules. Thus if $\rho: M \rightarrow \mathcal{H} \otimes M$, is the structure map of a left $\mathcal{H}$-comodule $M$, 
we write $\rho(m)=m^{(\bar{1})} \otimes m^{(\bar{0})}$, where
summation is understood. 

 A coalgebra $C$ is called a left $\mathcal{H}$-comodule coalgebra if $C$ is a left $\mathcal{H}$-comodule and the comultiplication map $ C \rightarrow C \otimes C $ and the 
counit map $C \rightarrow k$ are morphisms of $\mathcal{H}$
-comodules, i.e., for all $a \in C,$ 
\begin{eqnarray*}
& & \sum a^{(0) \bar{(1)}} a^{(1) \bar{(1)}} \otimes a^{(0) \bar{(0)}} \otimes a^{(1) \bar{(0)}}
= \sum a^{\bar{(1)}} \otimes a^{\bar{(0)} (0)} \otimes a^{\bar{(0)} (1)},  
\end{eqnarray*} 
 \[
\sum a^{\bar{(1)}} \epsilon(a^{(\bar{0})}) = \epsilon(a) 1_{\mathcal{H}}, \]
where $\epsilon$ is the counit of $C$. The reader can consult~\cite{cdr96} for various examples of comodule coalgebras.\\

If $\mathcal{H}$ is a bialgebra and $C$ a left $\mathcal{H}$-comodule coalgebra, the
crossed coproduct coalgebra $C \!> \! \blacktriangleleft \!\mathcal{H}$ is defined as $C \otimes \mathcal{H} $ with the coalgebra structure 
\begin{eqnarray*}
& &\Delta(a \otimes g) = \sum a^{(0)} \otimes a^{(1) \bar{(1)}} g^{(0)} \otimes a^{(1) \bar{(0)}} \otimes g^{(1)},\\
& &\epsilon (a \otimes g) = \epsilon(a) \epsilon(g),
\end{eqnarray*} 
for $g \in \mathcal{H}$ and $a \in C.$\\

For more information on cyclic homology and Hopf algebras we recommend
\cite{aC994, ll92, sM95, sw69, cw94}.


\section{The Cocylindrical Module $\mathbf {C } \mathbf{\natural^{\textbf{}}} \mathbf{\mathcal{H}}$}
Let $\mathcal{H}$ be a Hopf algebra with bijective antipode and $C$ a left $\mathcal{H}$-comodule coalgebra. In this section we 
associate a cocylindrical module $C \natural^{\small{}} \mathcal{H}$ to this data by
$C \natural^{\small{}} \mathcal{H}= \{ \mathcal{H}^{\otimes(p+1)}\otimes C ^{\otimes (q+1)} \}_{p,q \ge 0} $.
We define the
operators $\tau_{p,q} ,\partial_{p,q} ,\sigma_{p,q}$ and $\bar{\tau}_{p,q} ,\bar{\partial}_{p,q} ,\bar{\sigma}_{p,q}$ 
as follows:

\begin{eqnarray} \label{eq:ccopcy1} \notag
& &\tau_{p,q}( g_0 , \dots , g_p \mid a_0 , \dots , a_q) 
= (S^{-1}(a_0^{\bar{(1)}}) \cdot (g_0, \dots , g_p) \mid
  a_1, \dots ,a_{q},a_0^{\bar{(0)}}),   \notag
\end{eqnarray}
\begin{eqnarray} \notag
& &\sigma_{p,q}^i(g_0,\dots,g_p \mid a_0 , \dots , a_q)  \\ \notag
& &= (g_0,\dots,g_p \mid a_0,\dots ,a_{i}, a_{i+2},\dots , a_q) \epsilon (a_{i+1}) \;\;\; 0 \le i < q,  \\ \notag
& &\partial_{p,q}^i(g_0,\dots,g_p \mid a_0 , \dots , a_q) \\ 
& &= (g_0,\dots,g_p \mid a_0,\dots , a_i^{(0)},a_i^{(1)}, a_{i+1},\dots , a_q) \;\;\; 0 \le i \le q, \\ \notag
& &\partial_{p,q}^{q+1}(g_0,\dots,g_p \mid a_0 , \dots , a_q) \\ \notag
& &= (S^{-1}(a_0^{(0) \bar{(1)}}) \cdot (g_1,\dots,g_p )\mid a_0^{(1)}, a_0,\dots , a_q,a_0^{(0) \bar{(0)}}),  \notag
\end{eqnarray}

where \;\;$h \cdot (g_0, \dots, g_p) :=( h^{(0)}  g_0, \dots ,h^{(p)} g_p ),$ and we see that
\begin{eqnarray*}
& & S^{-1}(a_0^{\bar{(1)}}) \cdot (g_0, \dots , g_p)= (S^{-1}(a_0^{\bar{(1)}(p)}) g_0 ,\dots,S^{-1}(a_0^{\bar{(1)}(0)}) g_p) \\
&=& (S^{-1}(a_0^{\bar{(1)}}) g_0 ,\dots,S^{-1}(a_0^{(\overline{p+1})}) g_p). 
\end{eqnarray*}
\begin{eqnarray} \label{eq:ccopcy2} \notag
& &\bar{\tau}_{p,q}( g_0 , \dots , g_p \mid a_0 , \dots , a_q) \notag
=(g_1, \dots , g_{p-1},(a_0^{\bar{(1)}} \dots a_q^{\bar{(1)}})  g_0 \mid  a_0^{\bar{(0)}}, \dots
, a_{q}^{\bar{(0)}}), \\ \notag
& &\bar{\sigma}_{p,q}^i(g_0,\dots,g_p \mid a_0 , \dots , a_q)\\ \notag
& &= (g_0,\dots ,g_{i}, g_{i+2}, \dots,,g_{p} \mid a_0 , \dots, a_q) \epsilon(g_{i+1}) \;\; 0 \le i < p,\\
& &\bar{\partial}_{p,q}^i (g_0,\dots,g_p \mid a_0 , \dots , a_q)\\ \notag
& &=(g_0,\dots,g_i^{(0)},g_{i}^{(1)},\dots,g_p \mid a_0,\dots , a_q) \;\;\; 0 \le i \le p, \\ \notag
& &\bar{\partial}_{p,q}^{p+1} (g_0,\dots,g_p \mid a_0 , \dots , a_q)\\ \notag
& &=(g_{0}^{(1)},g_0,\dots,g_p,(a_0^{\bar{(1)}} \dots a_q^{\bar{(1)}}) g_{0}^{(0)}\mid a_0^{\bar{(0)}},\dots , a_q^{\bar{(0)}}). \notag
\end{eqnarray} 

\begin{theorem}
$C \natural^{\text{}} \mathcal{H}$ with the operators defined in (\ref{eq:ccopcy1}),(\ref{eq:ccopcy2}) is a cocylindrical module.  
\end{theorem}
\begin{proof}

We verify only some of the important identities. The other identities are similar to check.\\

\item[$\bullet \quad$]$  \tau_{p,q+1} \partial^0_{p,q} = \partial^{q+1}_{p,q}$\\

$ \tau_{p,q+1} \partial^0_{p,q} (g_0 , \dots , g_p \mid a_0 , \dots , a_q)$ 
\begin{eqnarray*}
&=&\tau_{p,q+1}(g_0, \dots , g_p \mid a_0^{(0)},a_0^{(1)}, \dots,a_{q} )\\
&=&\partial_{p,q} (S^{-1}(a_0^{(0) \bar{(1)}}) \cdot (g_0, \dots , g_p) \mid a_0^{(1)}, \dots,a_{q},a_0^{(0) \bar{(0)}} )\\
&=&\partial^{q+1}_{p,q} (g_0 , \dots , g_p \mid a_0 , \dots , a_q).
\end{eqnarray*}
\item[$\bullet \quad$]$ \tau_{p,q-1} \sigma^0_{p,q} =  \sigma^{q+1}_{p,q} {(\tau^{q}_{p,q})}^2 $
\begin{eqnarray*}
\tau_{p,q-1} \sigma^0_{p,q} (g_0 , \dots , g_p \mid a_0 , \dots , a_q)
=(S^{-1}(a_0^{\bar{(1)}}) \cdot (g_0,\dots,g_p) \mid a_1, \dots ,a_q,a_0^{\bar{(0)}}) \epsilon(a_1), 
\end{eqnarray*}
$\sigma^{q+1}_{p,q} {(\tau_{p,q})}^2  (g_0 , \dots , g_p \mid a_0 , \dots , a_q)$
\begin{eqnarray*}
&=&\sigma^{q+1}_{p,q} ((S^{-1}(a_1^{\bar{(1)})})(S^{-1}(a_0^{\bar{(1)})})) \cdot (g_0,\dots,g_p) \mid a_2, \dots ,a_q,a_0^{\bar{(0)}},a_1^{\bar{(0)}}) \\
&=&(S^{-1}(a_0^{\bar{(1)}}) \cdot (g_0,\dots,g_p) \mid a_1, \dots ,a_q,a_0^{\bar{(0)}}) \epsilon(a_1).\\ 
\end{eqnarray*}

\item[$\bullet \quad$]$ \bar{\tau}_{p-1,q} \bar{\sigma}^0_{p,q} = \bar{\sigma}^{p-1}_{p,q} (\bar{\tau}_{p,q})^2 $\\
\begin{eqnarray*}
\bar{\tau}_{p-1,q} \bar{\sigma}^0_{p,q} (g_0 , \dots , g_p \mid a_0 , \dots , a_q)=
(g_2,\dots,g_p,(a_0^{\bar{(1)}} \dots a_q^{\bar{(1)}})g_0 \mid a_0^{\bar{(0)}}, \dots ,a_q^{\bar{(0)}}) \epsilon(g_1),
\end{eqnarray*}
$ \bar{\sigma}^{p-1}_{p,q} (\bar{\tau}_{p,q})^2  (g_0 , \dots , g_p \mid a_0 , \dots , a_q)$ 
\begin{eqnarray*}
&=&\bar{\sigma}^{p-1}_{p,q} (g_2,\dots,g_p,(a_0^{\bar{(2)}} \dots a_q^{\bar{(2)}})g_0,
(a_0^{\bar{(1)}} \dots a_q^{\bar{(1)}})g_1 \mid a_0^{\bar{(0)}}, \dots ,a_q^{\bar{(0)}})\\
&=&(g_2,\dots,g_p,(a_0^{\bar{(1)}} \dots a_q^{\bar{(1)}})g_0 \mid a_0^{\bar{(0)}}, \dots ,a_q^{\bar{(0)}}) \epsilon(g_1).
\end{eqnarray*}
\item[$\bullet \quad$]$ \bar{\tau}_{p+1,q} \bar{\partial}^0_{p,q} = \bar{\partial}^{q+1}_{p,q}$\\
$ \bar{\tau}_{p+1,q} \bar{\partial}^0_{p,q} (g_0 , \dots , g_p \mid a_0 , \dots , a_q)$ 
\begin{eqnarray*}
&=& (g_0^{(1)}, g_1, \dots ,g_{p-1},( a_0^{\bar{(1)}} \dots a_q^{\bar{(1)}})  g_0^{(0)}
\mid a_0^{\bar{(0)}}, \dots,a_q^{\bar{(0)}})\\   
&=& \bar{\partial}_{q+1}^{p,q} (g_0 , \dots , g_p \mid a_0 , \dots , a_q).\\
\end{eqnarray*}
\item[$\bullet \quad$] $\tau_{p,q} \bar{\tau}_{p,q} = \bar{\tau}_{p,q} \tau_{p,q}$\\
$(\tau_{p,q} \bar{\tau}_{p,q})(g_0,\dots,g_p \mid a_0,\dots ,a_q)$
\begin{eqnarray*}
&=&\tau_{p,q}(  g_1,\dots,g_{p},(a_0^{\bar{(1)}} \dots a_q^{\bar{(1)}}) g_0 
\mid a_0^{\bar{(0)}}, \dots,a_q^{\bar{(0)}})\\
&=&(S^{-1}(a_0^{\bar{(1)}}) \cdot ( g_1, 
\dots g_{p},(a_0^{\bar{(2)}} \dots a_q^{\bar{(1)}})  g_0) \mid a_1^{\bar{(0)}},  \dots a_{q}^{\bar{(0)}},a_0^{\bar{(0)}} )\\
&=&(
S^{-1}(a_0^{\bar{(1)} (0)}) \cdot (g_1, 
\dots g_{p}),S^{-1}(a_0^{\bar{(1)} (1)})  (a_0^{\bar{(1)} (2)}a_{1}^{\bar{(1)}} \dots a_q^{\bar{(1)}})  g_0) \mid a_1^{\bar{(0)}}, \dots a_{q}^{\bar{(0)}},a_0^{\bar{(0)}} )\\
&=&(
S^{-1}(a_0^{\bar{(1)} }) \cdot (g_1, 
\dots g_{p}),(a_{1}^{\bar{(1)}} \dots a_q^{\bar{(1)}}) g_0 \mid a_1^{\bar{(0)}}, \dots a_{q}^{\bar{(0)}}, a_0^{\bar{(0)}} ).\\
\end{eqnarray*} 
$(\bar{\tau}_{p,q} \tau_{p,q})(g_0,\dots,g_p \mid a_0,\dots ,a_q)$
\begin{eqnarray*}
&=& \bar{\tau}_{p,q}(S^{-1}(a_0^{\bar{(1)}}) \cdot (g_0, \dots , g_p) \mid a_1, \dots,a_{q},a_0^{\bar{(0)}})\\
&=& \bar{\tau}_{p,q}(S^{-1}(a_0^{\bar{(1)}(0)})  g_0, S^{-1}(a_0^{\bar{(1)}(1)}) \cdot (g_1,\dots , g_{p}) \mid
  a_1, \dots,a_{q},a_0^{\bar{(0)}})\\
&=& ( S^{-1}(a_0^{\bar{(2)}(1)}) \cdot (g_1,\dots , g_{p}),(a_1^{\bar{(1)}} \dots a_p^{\bar{(1)}} a_0^{\bar{(1)}}) S^{-1}(a_0^{\bar{(2)}(1)})  g_0) \mid
  a_1, \dots,a_{q},a_0^{\bar{(0)}})\\
&=& ( S^{-1}(a_0^{\bar{(1)}(1)}) \cdot (g_1,\dots , g_{p}),(a_1^{\bar{(1)}} \dots a_p^{\bar{(1)}} a_0^{\bar{(1)}(0)})  S^{-1}(a_0^{\bar{(1)}(1)}) g_0) \mid
  a_1, \dots,a_{q},a_0^{\bar{(0)}})\\
&=& ( S^{-1}(a_0^{\bar{(1)}}) \cdot (g_1,\dots , g_{p}),(a_1^{\bar{(1)}} \dots a_p^{\bar{(1)}})  g_0) \mid
  a_1, \dots,a_{q},a_0^{\bar{(0)}}).\\
\end{eqnarray*}
\item[$\bullet \quad$] $\bar{\tau}_{p,q}^{p+1} \;
\tau_{p,q}^{q+1}=\tau_{p,q}^{q+1} \;
\bar{\tau}_{p,q}^{p+1}=id_{p,q}$\\
$\tau_{p,q}^{q+1} \bar{\tau}_{p,q}^{p+1} (g_0,\dots,g_p \mid a_0, \dots , a_q) $
\begin{eqnarray*}
&=&\tau^{q+1}_{p,q} ( (a_0^{(\overline{p+1})} \dots a_q^{(\overline{p+1})} )  g_0 , \dots ,
(a_0^{\bar{(1)}} \dots a_q^{\bar{(1)}} )  g_p \mid a_0^{\bar{(0)}} , \dots , a_q^{\bar{(0)}} )\\
&=&S^{-1}(a_0^{(\overline{1})} \dots a_q^{(\overline{1})}) \cdot ((a_0^{(\overline{p+2})} \dots a_q^{(\overline{p+2})} )  g_0 , \dots ,
(a_0^{\bar{(2)}} \dots a_q^{\bar{(2)}} )  g_p) \mid a_0^{\bar{(0)}} , \dots , a_q^{\bar{(0)}} )\\
&=&(g_0,\dots,g_p \mid a_0, \dots , a_q).
\end{eqnarray*}
\end{proof}
\begin{corollary} The diagonal $\Delta (C \natural^{} \mathcal{H})$ is a cocyclic module.
\end{corollary}

\section{Relation of $\mathbf{\Delta (C \natural^{} \mathcal{H})}$ with the Cocyclic Module of the 
Crossed Coproduct Coalgebra $\mathbf{ C \!> \! \blacktriangleleft \!\mathcal{H} }$} 
We define a map $\phi : \mathsf{C}^{\bullet}(C \!> \! \blacktriangleleft \!\mathcal{H}) \rightarrow \Delta(C \natural^{} \mathcal{H})$, by 
$$\phi ( a_0 \otimes g_0,\dots ,a_n \otimes g_n )=$$
\begin{eqnarray*}
(S^{-1}(a_1^{\bar{(1)}} \dots a_p^{\bar{(1)}})  g_0 ,S^{-1}(a_2^{\bar{(2)}} \dots a_p^{\bar{(2)}}) g_1, \dots ,S^{-1}(a_{p}^{(\overline{p})})  g_{p-1},g_p \mid a_0,a_1^{\bar{(0)}} , \dots , a_q^{\bar{(0)}}),  
\end{eqnarray*}
where 
$\mathsf{C}^{\bullet}(C  \!> \! \blacktriangleleft \! \mathcal{H})$ denotes the cocyclic module of the crossed coproduct coalgebra $ C  \!> \! \blacktriangleleft  \!\mathcal{H} $.\\
\begin{theorem} $\phi:\mathsf{C}^{\bullet}(C \!> \! \blacktriangleleft \! \mathcal{H}) \rightarrow \Delta (C \natural^{} \mathcal{H})$ defines a cocyclic map.
\end{theorem}
\begin{proof}
We show that $\phi$ commutes with the cosimplicial and cyclic operators:
$\text{}\\$
\item[$\bullet \quad$] $  \tau_{n,n} \bar{\tau}_{n,n} \phi = \phi \tau_n^{C > \! \blacktriangleleft \mathcal{H}} $ 
\begin{eqnarray*}
& & \tau_{n,n}\bar{\tau}_{n,n} \phi (a_0 \otimes g_0, \dots , a_n \otimes g_n)\\
&=& \tau_{n,n} \bar{\tau}_{n,n}  (S^{-1}(a_1^{\bar{(1)}} \dots a_n^{\bar{(1)}}) g_0 ,
S^{-1}(a_2^{\bar{(2)}} \dots a_n^{\bar{(2)}})  g_1, \dots ,S^{-1}(a_n^{(\overline{n})}) 
 g_{n-1},g_n \mid a_0,a_1^{\bar{(0)}} , \dots , a_n^{\bar{(0)}})\\ 
&=& \tau_{n,n} (
S^{-1}(a_2^{\bar{(3)}} \dots a_n^{\bar{(3)}})  g_1, \dots ,S^{-1}(a_n^{(\overline{n+2})}) g_n ,(a_0^{\bar{(1)}} \dots a_n^{\bar{(1)}})
 S^{-1}(a_1^{\bar{(2)}} \dots a_n^{\bar{(2)}})  g_0 \\  
& & \hspace{11cm} \mid a_0^{\bar{(0)}} ,a_1, \dots , a_n^{\bar{(0)}})\\
&=& \tau_{n,n} (
S^{-1}(a_2^{\bar{(1)}} \dots a_n^{\bar{(1)}})  g_1, \dots , g_n ,a_0^{\bar{(1)}}g_0   
\mid a_0^{\bar{(0)}}, \dots , a_n^{\bar{(0)}})\\ 
&=& (S^{-1}(a_0^{\bar{(1)}}) \cdot(
S^{-1}(a_2^{\bar{(1)}} \dots a_n^{\bar{(1)}})  g_1, \dots ,g_n ,a_0^{\bar{(2)}}g_0   )   
\mid a_1,a_2^{\bar{(0)}}, \dots , a_n^{\bar{(0)}},a_0^{\bar{(0)}} )\\ 
&=& (
S^{-1}(a_2^{\bar{(1)}} \dots a_n^{\bar{(1)}} a_0^{\bar{(1)}})  g_1, \dots ,
S^{-1}(a_0^{(\overline{n})}) g_n ,S^{-1}(a_0^{(\overline{n+1})}) a_0^{(\overline{n+2})} g_0)   
\mid a_1,a_2^{\bar{(0)}}, \dots , a_n^{\bar{(0)}},a_0^{\bar{(0)}} )\\
&=& (
S^{-1}(a_2^{\bar{(1)}} \dots a_n^{\bar{(1)}} a_0^{\bar{(1)}})  g_1, \dots ,
S^{-1}(a_0^{(\overline{n})}) g_n , g_0)   
\mid a_1,a_2^{\bar{(0)}}, \dots , a_n^{\bar{(0)}},a_0^{\bar{(0)}} )\\
&=&\phi (a_1 \otimes g_1, \dots , a_n \otimes g_n, a_0 \otimes g_0)=
\phi \tau_n^{C \!> \! \blacktriangleleft \!\mathcal{H}} ( a_0 \otimes g_0,\dots ,a_n \otimes g_n ).     
\end{eqnarray*}
\item[$\bullet \quad$] $ \sigma^i_{n,n} \bar{\sigma}^i_{n,n}  \phi = \phi \sigma_i^{C > \! \blacktriangleleft \mathcal{H}}\;\;\; 0 \le i < n $ 
\begin{eqnarray*}
& & \sigma^i_{n,n} \bar{\sigma}^i_{n,n}   \phi (a_0 \otimes g_0 , \dots , a_n \otimes g_{n}) \\
&=& \sigma^i_{n,n} \bar{\sigma}^i_{n,n}   (S^{-1}(a_1^{\bar{(1)}} \dots a_n^{\bar{(1)}})  g_0 ,
S^{-1}(a_2^{\bar{(2)}} \dots a_n^{\bar{(2)}})  g_1, \dots ,S^{-1}(a_n^{(\overline{n})}) 
 g_{n-1},g_n \mid a_0,a_1^{\bar{(0)}} , \dots , a_n^{\bar{(0)}})\\
&=& \sigma^i_{n,n} (S^{-1}(a_1^{\bar{(1)}} \dots a_n^{\bar{(1)}})  g_0 ,
S^{-1}(a_2^{\bar{(2)}} \dots a_n^{\bar{(2)}})  g_1, \dots ,S^{-1}(a_n^{(\overline{n})}) 
 g_{n-1},g_n \\
& & \hspace{6.5cm} \mid a_0,a_1^{\bar{(0)}} , \dots , a_n^{\bar{(0)}}) \epsilon(S^{-1}(a_{i+2}^{(\overline{i+2})} \dots a_n^{(\overline{i+2})}) g_{i+1}) \\
&=& (S^{-1}(a_1^{\bar{(1)}} \dots a_{i-1}^{\bar{(1)}}  a_{i+1}^{\bar{(1)}} \dots  a_n^{\bar{(1)}})  g_0 ,
S^{-1}(a_2^{\bar{(2)}} \dots a_{i-1}^{\bar{(2)}}  a_{i+1}^{\bar{(2)}} \dots a_n^{\bar{(2)}})  g_1, \dots ,S^{-1}(a_n^{(\overline{n})}) 
 g_{n-1},g_n \\
& & \hspace{6.5cm} \mid a_0,a_1^{\bar{(0)}} , \dots ,a_i^{\bar{(0)}},a_{i+2}^{\bar{(0)}},\dots, a_n^{\bar{(0)}}) \epsilon(g_{i+1}) \epsilon(a_{i+1})\\
&=& (S^{-1}(a_1^{\bar{(1)}} \dots a_{i-1}^{\bar{(1)}}  a_{i+1}^{\bar{(1)}} \dots  a_n^{\bar{(1)}})  g_0 ,
S^{-1}(a_2^{\bar{(2)}} \dots a_{i-1}^{\bar{(2)}}  a_{i+1}^{\bar{(2)}} \dots a_n^{\bar{(2)}})  g_1, \dots ,S^{-1}(a_n^{(\overline{n})}) 
 g_{n-1},g_n \\
& & \hspace{6.5cm} \mid a_0,a_1^{\bar{(0)}} , \dots ,a_i^{\bar{(0)}},a_{i+2}^{\bar{(0)}},\dots, a_n^{\bar{(0)}}) \epsilon(g_{i+1} \otimes a_{i+1})\\
&=& \phi \sigma_i^{C \!> \! \blacktriangleleft \!\mathcal{H}} (a_0 \otimes g_0 , \dots , a_n \otimes g_n).\\ 
\end{eqnarray*}

\item[$\bullet \quad$] $ \partial^{n+1}_{n,n} \bar{\partial}^{n+1}_{n,n} \phi = \phi \partial_{n+1}^{C \!> \! \blacktriangleleft \!\mathcal{H}} $ 
\begin{eqnarray*}
& & \partial^{n+1}_{n,n} \bar{\partial}^{n+1}_{n,n}   \phi (a_0 \otimes g_0 , \dots , a_n \otimes g_n) \\
&=& \partial^{n+1}_{n,n} \bar{\partial}^{n+1}_{n,n}   (S^{-1}(a_1^{\bar{(1)}} \dots a_n^{\bar{(1)}})  g_0 ,
S^{-1}(a_2^{\bar{(2)}} \dots a_n^{\bar{(2)}})  g_1, \dots ,S^{-1}(a_n^{(\overline{n})}) 
 g_{n-1},g_n \\
& &\hspace{10cm} \mid a_0,a_1^{\bar{(0)}} , \dots , a_n^{\bar{(0)}})
\end{eqnarray*}
\begin{eqnarray*}
&=& \partial^{n+1}_{n,n}  (S^{-1}(a_1^{\bar{(2)}(0)} \dots a_n^{\bar{(2)}(0)})  g_0^{(1)} ,
S^{-1}(a_2^{\bar{(3)}} \dots a_n^{\bar{(3)}})  g_1, \dots , 
 g_n,  \\
& & \hspace{4cm} (a_0^{\bar{(1)}} \dots a_n^{\bar{(1)}}) S^{-1}(a_1^{\bar{(2)}(1)} \dots a_n^{\bar{(2)}(1)})  g_0^{(0)} 
\mid a_0^{\bar{(0)}} , \dots , a_n^{\bar{(0)}})\\
&=& \partial^{n+1}_{n,n}  (S^{-1}(a_0^{\bar{(3)}} \dots a_n^{\bar{(3)}})  g_0^{(1)} ,
S^{-1}(a_1^{\bar{(4)}} \dots a_n^{\bar{(4)}})  g_1, \dots , 
 g_n,  \\
& & \hspace{4.5cm} (a_0^{\bar{(1)}} \dots a_n^{\bar{(1)}}) S^{-1}(a_2^{\bar{(2)}} \dots a_n^{\bar{(2)}})  g_0^{(0)} 
\mid a_0^{\bar{(0)}} , \dots , a_n^{\bar{(0)}})\\
&=& \partial^{n+1}_{n,n}  (S^{-1}(a_1^{\bar{(1)}} \dots a_n^{\bar{(1)}})  g_0^{(1)} ,
S^{-1}(a_2^{\bar{(2)}} \dots a_n^{\bar{(2)}})  g_1, \dots , 
 g_n, a_0^{\bar{(1)}} g_0^{(0)} 
\mid a_0^{\bar{(0)}} , \dots , a_n^{\bar{(0)}})
\end{eqnarray*}
\begin{eqnarray*}
&=& (S^{-1}(a_0^{(0) \bar{(1)}}) \cdot (S^{-1}(a_1^{\bar{(1)}} \dots a_n^{\bar{(1)}})  g_0^{(1)} ,
S^{-1}(a_2^{\bar{(2)}} \dots a_n^{\bar{(2)}})  g_1, \dots , 
 g_n, a_0^{(0) \bar{(2)}} a_0^{(1) \bar{(1)}} g_0^{(0)}) \\
& & \hspace{8.5cm} \mid a_0^{(1) \bar{(0)}} , \dots , a_n^{\bar{(0)}},a_0^{(0) \bar{(0)}} )\\
&=& ( S^{-1}(a_1^{\bar{(1)}} \dots a_n^{\bar{(1)}} a_0^{(0) \bar{(1)}})  g_0^{(0)} ,
S^{-1}(a_2^{\bar{(2)}} \dots a_n^{\bar{(2)}}a_0^{(0) \bar{(2)}})  g_1, \dots , 
 S^{-1}(a_0^{(0) (\overline{n+1})}) g_n,\\ 
& & \hspace{3.5cm} S^{-1}(a_0^{(0) (\overline{n+2})}) a_0^{(0) (\overline{n+3})} a_0^{(1) \bar{(1)}} g_0^{(0)}) 
\mid a_0^{(1) \bar{(0)}} , \dots , a_n^{\bar{(0)}},a_0^{(0) \bar{(0)}} )\\
&=& ( S^{-1}(a_1^{\bar{(1)}} \dots a_n^{\bar{(1)}} a_0^{(0) \bar{(1)}})  g_0^{(1)} ,
S^{-1}(a_2^{\bar{(2)}} \dots a_n^{\bar{(2)}}a_0^{(0) \bar{(2)}})  g_1, \dots , 
 S^{-1}(a_0^{(0) (\overline{n+1})}) g_n,\\ 
& & \hspace{6.5cm}  a_0^{(1) \bar{(1)}} g_0^{(0)}) 
\mid a_0^{(1) \bar{(0)}} , \dots , a_n^{\bar{(0)}},a_0^{(0) \bar{(0)}} ).
\end{eqnarray*}

$ \phi \partial_{n+1}^{C > \! \blacktriangleleft  \mathcal{H}} (a_0 \otimes g_0 , \dots , a_n \otimes g_n)$ 
\begin{eqnarray*}
&=&\phi ( a_0^{(1) \bar{(0)}} \otimes g_0^{(1)},a_1 \otimes g_1,\dots,a_n \otimes g_n,
 a_0^{(0)} \otimes a_0^{(1) \bar{(1)}} g_0^{(0)})\\
&=& (S^{-1}( a_1^{\bar{(1)}} \dots a_n^{\bar{(1)}} a_0^{(0) \bar{(1)}} ) g_0^{(1)} ,
S^{-1}(a_2^{\bar{(2)}} \dots  a_n^{\bar{(2)}} a_0^{(0) \bar{(2)}} )  g_1, \dots , \\
& & \hspace{2cm} S^{-1}(a_0^{(0) (\overline{n+1})}) g_n ,a_0^{(1) \bar{(1)}} g_0^{(0)}) \mid a_0^{(1) \bar{(0)}} , \dots , a_n^{\bar{(0)}},a_0^{(0) \bar{(0)}} )\\ 
\end{eqnarray*}

\item[$\bullet \quad$] $ \partial^i_{n,n} \bar{\partial}^i_{n,n} \phi = \phi \partial_i^{C > \! \blacktriangleleft \mathcal{H}}\;\;\; 0 \le i \le n $ 
\begin{eqnarray*}
& & \partial^i_{n,n} \bar{\partial}^i_{n,n}   \phi (a_0 \otimes g_0 , \dots , a_n \otimes g_n) \\
&=& \partial^i_{n,n} \bar{\partial}^i_{n,n}   (S^{-1}(a_1^{\bar{(1)}} \dots a_n^{\bar{(1)}})  g_0 ,
S^{-1}(a_2^{\bar{(2)}} \dots a_n^{\bar{(2)}})  g_1, \dots ,S^{-1}(a_n^{(\overline{n})}) 
 g_{n-1},g_n  \mid a_0,a_1^{\bar{(0)}} , \dots , a_n^{\bar{(0)}})\\
&=& \partial^i_{n,n} (S^{-1}(a_1^{\bar{(1)}} \dots a_n^{\bar{(1)}})  g_0 ,
S^{-1}(a_2^{\bar{(2)}} \dots a_n^{\bar{(2)}})  g_1, \dots ,S^{-1}(a_{i+1}^{(\overline{i+1})(1)} \dots a_{i+1}^{(\overline{i+1})(1)}) g_i^{(0)},\\ 
& & \hspace{3cm} S^{-1}(a_{i+1}^{(\overline{i+1})(0)} \dots a_{i+1}^{(\overline{i+1})(0)}) g_i^{(1)},\dots,
S^{-1}(a_n^{(\overline{n})}) 
 g_{n-1},g_n \mid a_0,a_1^{\bar{(0)}} , \dots , a_n^{\bar{(0)}})\\
&=& \partial^i_{n,n} (S^{-1}(a_1^{\bar{(1)}} \dots a_n^{\bar{(1)}})  g_0 ,
S^{-1}(a_2^{\bar{(2)}} \dots a_n^{\bar{(2)}})  g_1, \dots ,S^{-1}(a_{i+1}^{(\overline{i+1})} \dots a_n^{(\overline{i+1})}) g_i^{(0)},\\ 
& & \hspace{3cm} S^{-1}(a_{i+1}^{(\overline{i+2})} \dots a_n^{(\overline{i+2})}) g_i^{(1)},\dots,
S^{-1}(a_n^{(\overline{n+1})}) 
 g_{n-1},g_n \mid a_0,a_1^{\bar{(0)}} , \dots , a_n^{\bar{(0)}})\\
&=& (S^{-1}(a_1^{\bar{(1)}} \dots a_n^{\bar{(1)}})  g_0 ,
S^{-1}(a_2^{\bar{(2)}} \dots a_n^{\bar{(2)}})  g_1, \dots ,S^{-1}(a_{i+1}^{(\overline{i+1})} \dots a_n^{(\overline{i+1})}) g_i^{(0)},\\ 
& & \hspace{.2cm} S^{-1}(a_{i+1}^{(\overline{i+2})} \dots a_n^{(\overline{i+2})}) g_i^{(1)},\dots,
S^{-1}(a_n^{(\overline{n+1})}) 
 g_{n-1},g_n \mid a_0,a_1^{\bar{(0)}} , \dots , a_i^{\bar{(0)}(0)},a_i^{\bar{(0)}(1)}, \dots, a_n^{\bar{(0)}}).
\end{eqnarray*}
$ \phi \partial_i^{C > \! \blacktriangleleft \mathcal{H}} (a_0 \otimes g_0 , \dots , a_n \otimes g_n)$ 
\begin{eqnarray*}
&=&\phi (a_0 \otimes g_0 , \dots ,a_i^{(0)} \otimes a_i^{(1) \bar{(1)}} g_i^{(0)}, a_i^{(1) \bar{(0)}} 
\otimes g_i^{(1)},\dots ,a_n \otimes g_n)\\
&=& (S^{-1}(a_1^{\bar{(1)}} \dots a_i^{(0) \bar{(1)}}a_i^{(1) \bar{(1)}} \dots a_n^{\bar{(1)}})  g_0 ,
S^{-1}(a_2^{\bar{(2)}} \dots a_i^{(0) \bar{(2)}}a_i^{(1) \bar{(2)}} \dots a_n^{\bar{(2)}})  g_1, \dots ,\\
& & \hspace{.1cm} S^{-1}(a_i^{(0) (\overline{i})}a_i^{(1) (\overline{i})}  \dots a_n^{(\overline{i})})g_{i-1}
,S^{-1}(a_i^{(1) (\overline{i+1})}  \dots a_n^{(\overline{i+1})}) a_i^{(1) (\overline{i+2})} g_i^{(0)} ,\\
& & \hspace{.4cm} S^{-1}(a_{i+1}^{(\overline{i+2})}  \dots a_n^{(\overline{i+2})})  g_i^{(1)} ,\dots, S^{-1}(a_n^{(\overline{n+1})}) 
 g_{n-1},g_n \mid a_0,a_1^{\bar{(0)}} , \dots ,a_i^{(0) \bar{(0)}},a_i^{(1) \bar{(0)}},\dots a_n^{\bar{(0)}})\\ 
&=& (S^{-1}(a_1^{\bar{(1)}} \dots a_i^{ \bar{(1)}} \dots a_n^{\bar{(1)}})  g_0 ,
S^{-1}(a_2^{\bar{(2)}} \dots a_i^{ \bar{(2)}} \dots a_n^{\bar{(2)}})  g_1, \dots ,S^{-1}(a_i^{(\overline{i})}  \dots a_n^{(\overline{i})})g_{i-1}
\\
& & \hspace{1cm} ,S^{-1}(a_{i+1}^{(\overline{i+1})}  \dots a_n^{(\overline{i+1})}) g_i^{(0)} ,S^{-1}(a_{i+1}^{(\overline{i+2})}  
\dots a_n^{(\overline{i+2})})  g_i^{(1)} ,\dots, S^{-1}(a_n^{(\overline{n+1})}) 
 g_{n-1},g_n \\
& & \hspace{7.5cm} \mid a_0,a_1^{\bar{(0)}} , \dots ,a_i^{ \bar{(0)}(0)},a_i^{ \bar{(0)}(1)},\dots a_n^{\bar{(0)}}).\\ 
\end{eqnarray*}
\end{proof}
\begin{theorem} We have an isomorphism of cocyclic modules 
\[\Delta (C \natural^{} \mathcal{H}) \cong \mathsf{C}^{\bullet}(C  \!> \! \blacktriangleleft \!\mathcal{H}).\]
\end{theorem}
\begin{proof}
We define a map $\psi : \Delta(C \natural^{} \mathcal{H})  \rightarrow \mathsf{C}^{\bullet}(C \!> \! \blacktriangleleft \!\mathcal{H})$, by 
\begin{eqnarray*}
& &\psi (g_0, \dots , g_n \mid a_0 , \dots ,a_n)=\\
& &(a_0 \otimes ( a_1^{\bar{(1)}} \dots a_p^{\bar{(p)}}) g_0,a_1^{\bar{(0)}} \otimes ( a_2^{\bar{(2)}} \dots a_p^{(\overline{p-1})}) g_1,
\dots, a_{p-1}^{\bar{(0)}} \otimes ( a_p^{\bar{(0)}} ) g_{p-1}, a_p^{\bar{(0)}} \otimes g_p). 
\end{eqnarray*}
One can check that  $\psi$ is a cocyclic map and $\phi \circ \psi = \psi \circ \phi =id$. For example\\
$\phi \circ \psi ( g_0, \dots ,g_n \mid a_0, \dots , a_n)=$
\begin{eqnarray*}
& &(S^{-1}(a_1^{\bar{(1)}} \dots a_p^{\bar{(1)}})(a_1^{\bar{(2)}} a_2^{\bar{(4)}} \dots a_p^{(\overline{2p})})g_0,
S^{-1}(a_2^{\bar{(2)}} \dots a_p^{\bar{(2)}})(a_2^{\bar{(3)}} a_3^{\bar{(5)}} \dots a_p^{(\overline{2p-1})})g_1,\\ 
& & \hspace{5cm} \dots ,S^{-1}(a_p^{(\overline{p})})a_p^{(\overline{p+1})} g_{p-1},g_p \mid a_0,a_1^{\bar{(0)}},\dots,a_p^{\bar{(0)}})\\
& & =( g_0, \dots ,g_n \mid a_0, \dots , a_n).
\end{eqnarray*}
\end{proof}

Now, we can use the Eilenberg-Zilber theorem for cocylindrical modules to conclude $$H^{\bullet}(Tot(C \natural \mathcal{H}); \mathsf{W}) \cong 
HC^{\bullet} (\Delta( C \natural^{}  \mathcal{H}); \mathsf{W})\cong HC^{\bullet}(C  \!> \! \blacktriangleleft \!\mathcal{H}; \mathsf{W}).$$
 
\section{The Cocyclic Module $\mathbf{\mathsf{C}^{\bullet}_{\mathcal{H}^{}}}(C)$ of Coinvariants and Cyclic 
Cohomology of $C  \!> \! \blacktriangleleft \!\mathcal{H}$}

Let $M$ be a left $\mathcal{H}$-comodule with the coaction $\Delta_{\tiny{M}}$. We define 
the space of $p$-cochains on $\mathcal{H}$ with values in $M$ by $\mathsf{C}^p(\mathcal{H},M)=\mathcal{H}^{\otimes p} \otimes M$, and we define the coboundary $\boldsymbol{\delta} : \mathsf{C}^p(\mathcal{H},M) \rightarrow
\mathsf{C}^{p+1}(\mathcal{H},M)$ by\\ 

\begin{eqnarray}
& &\boldsymbol{\delta} (g_1,g_2,\dots , g_p,m)= 
(1,g_1,\dots ,g_p,m) +\\
& &\sum_{i=1}^{p} (-1)^i (g_1 , \dots ,g_i^{(0)}, g_{i}^{(1)},\dots , g_p ,m) + (-1)^{p+1} (g_1,\dots,g_{p},\Delta_{\tiny{M}} (m)). \notag
\end{eqnarray}
We denote the $p$-th cohomology of the complex $(\mathsf{C}^{\bullet}(\mathcal{H},M),\boldsymbol{\delta})$ by $H^p(\mathcal{H},M)$.\\

Now we define a left $\mathcal{H}$-coaction  on the first row of $C \natural^{} \mathcal{H}$, $C^{\natural^ {}}_{\mathcal{H}}$ 
= $ \{ \mathcal{H} \otimes C^{\otimes (n+1)} \}_{n \ge 0}$, by 
\begin{equation} \label{eq:ac}
\boldsymbol{\Delta} (g \mid a_0, \dots , a_n ) = ( S^{-1}(g^{(2)})(a_0^{\bar{(1)}} \dots a_q^{\bar{(1)}}) g^{(0)} \mid g^{(1)} \mid a_0^{\bar{(0)}}, \dots ,a_q^{\bar{(0)}}). 
\end{equation}
We define $\mathsf{C}^q(C_{\mathcal{H}}^{\natural^{}})= \mathcal{H}^{\otimes q} \otimes C_{\mathcal{H}}^{\natural^{}}$ and
a coaction on it by $ \boldsymbol{\Delta} (g_1, \dots , g_p \mid m) = (g_1, \dots , g_p \mid \boldsymbol{\Delta}(m))$ where
$m \in C^{\natural^{ }}_{\mathcal{H}}.$ So we can construct
$H^p(\mathcal{H},\mathsf{C}^q(C^{\natural}_{\mathcal{H}}))$.        
We define $\mathsf{C}^{\bullet}_{\mathcal{H}^{}}(C)$ as the coinvariant subspace of $C^{\natural^{}}_{\mathcal{H}}$ under the above
coaction, so that $\mathsf{C}^{n}_{\mathcal{H}^{}}(C)$ is the space of all $(g \mid a_0, \dots , a_n ) $ such that \[
\boldsymbol{\Delta} (g \mid a_0, \dots , a_n ) = (1 \mid g \mid a_0, \dots , a_n ). \]

Now we replace the complex $(C \natural^{} \mathcal{H},(\partial,\sigma,\tau),(\bar{\partial},\bar{\sigma},\bar{\tau}))$ with the complex 
$(\mathsf{C}^p(\mathcal{H},\mathsf{C}^q ( C^{\natural^{}}_{\mathcal{H}}),( \mathfrak{d},\mathfrak{s},\mathfrak{t}),(\bar{\mathfrak{d}},\bar{\mathfrak{s}},\bar{\mathfrak{t}}))$,
under the transformation defined by 
the maps $\boldsymbol{\beta} : (C \natural^{} \mathcal{H})_{p,q} \rightarrow \mathsf{C}^p(\mathcal{H},
\mathsf{C}^q(C_{\mathcal{H}}^{\natural^{}})),$
and $\boldsymbol{\gamma} :  \mathsf{C}^p(\mathcal{H},\mathsf{C}^q(C_{\mathcal{H}}^{\natural^{}})) \rightarrow  
(C \natural^{} \mathcal{H})_{p,q},$ defined by
\begin{eqnarray*}
& &\boldsymbol{\beta} ( g_0, \dots , g_p \mid a_0 , \dots , a_q ) = ( S^{-1}(g^{(1)}) \cdot (g_1, \dots , g_p )\mid g_0^{(0)} \mid
a_0, \dots , a_q ), \\
& &\boldsymbol{\gamma} ( g_1 , \dots , g_p \mid g \mid a_0 , \dots , a_q) = (g^{(0)},g^{(1)} \cdot (g_1 , \dots , g_p)  \mid 
a_0, \dots , a_q ).
\end{eqnarray*}

One can check that $\boldsymbol{\beta} \circ \boldsymbol{\gamma} = \boldsymbol{\gamma} \circ \boldsymbol{\beta} = id $. We find the operators
$( \mathfrak{d},\mathfrak{s},\mathfrak{t})$ and $(\bar{\mathfrak{d}},\bar{\mathfrak{s}},\bar{\mathfrak{t}})$
under this transformation. \\

$ \hspace{0.2cm} \bullet \quad $ 
First we compute $\bar{\mathfrak{b}} = \boldsymbol{\beta} \bar{b} \boldsymbol{\gamma}$. Since 
\begin{eqnarray*}
\bar{b} ( g_0, \dots , g_p \mid a_0 , \dots , a_q ) &=& \sum_{0 \le i \le p } (-1)^i (g_0 , \dots , g_i^{(0)}, g_i^{(1)} ,\dots ,
 g_p \mid a_0 , \dots , a_q ) \\
&+& (-1)^{(p+1)} (g_0^{(1)}, \dots , g_{p}, (a_q^{\bar{(1)}} \dots a_0^{\bar{(1)}}) g_0^{(0)} ) \mid  a_0^{\bar{(0)}}, \dots
, a_{q}^{\bar{(0)}}), 
\end{eqnarray*}
we see that,
\begin{eqnarray*}
& &\bar{\mathfrak{b}}( g_1 , \dots , g_p \mid g \mid a_0 , \dots , a_q) = \boldsymbol{\beta} \bar{b} (g^{(0)},g^{(1)} \cdot
(g_1, \dots ,g_p) \mid a_0 , \dots , a_q )\\
&=& \boldsymbol{\beta} \{ (g^{(0)},g^{(1)}, g^{(2)} \cdot (g_1, \dots ,g_p) \mid a_0 , \dots , a_q )\\ 
&+& \sum _{1 \le i \le p} (-1)^i ( g^{(0)},g^{(1)}g_1,\dots,g^{(i)}g_i^{(0)},g^{(i+1)}g_i^{(1)},\dots,g^{(p+1)}g_p\mid a_0 , \dots , a_q)\\ 
&+& (-1)^{(p+1)} (g^{(1)},g^{(2)} \cdot (g_1, \dots ,g_p),(a_0^{\bar{(1)}} \dots a_q^{\bar{(1)}})g^{(0)}  \mid a_0^{\bar{(0)}} , \dots , a_q^{\bar{(0)}} )\\
&=& ( S^{-1}(g^{(1)}) \cdot (g^{(2)}, g^{(3)} \cdot (g_1, \dots ,g_p)) \mid g^{(0)} \mid a_0 , \dots , a_q )\\ 
&+& \sum _{1 \le i \le p} (-1)^i (S^{-1}(g^{(1)}) \cdot ( g^{(2)}g_1,\dots,g^{(i+1)}g_i^{(0)},g^{(i+2)}g_i^{(1)},\dots,g^{(p+2)}g_p )\mid g^{(0)} \mid a_0 , \dots , a_q)\\ 
&+& (-1)^{(p+1)} (S^{-1}(g^{(2)}) \cdot (g^{(3)} \cdot (g_1, \dots ,g_p),(a_0^{\bar{(1)}} \dots a_q^{\bar{(1)}})g^{(0)})\mid g^{(1)} \mid a_0^{\bar{(0)}} , \dots , a_q^{\bar{(0)}} )\\
&=& ( S^{-1}(g^{(2)}) g^{(3)},(S^{-1}(g^{(1)}) g^{(4)}) \cdot (g_1, \dots ,g_p) \mid g^{(0)} \mid a_0 , \dots , a_q )\\ 
\end{eqnarray*}
\begin{eqnarray*}
&+& \sum _{1 \le i \le p} (-1)^i (S^{-1}(g^{(1)}) \cdot ( g^{(2)} \cdot (g_1,\dots, g_i^{(0)}, g_i^{(1)},\dots, g_p ))\mid g^{(0)} \mid a_0 , \dots , a_q)\\ 
&+& (-1)^{(p+1)} ( (S^{-1}(g^{(3)}) (g^{(4)}) \cdot (g_1, \dots ,g_p),S^{-1}(g^{(2)})(a_0^{\bar{(1)}} \dots a_q^{\bar{(1)}})g^{(0)})\mid g^{(1)} \mid a_0^{\bar{(0)}} , \dots , a_q^{\bar{(0)}} )\\
&=& ( 1 , g_1, \dots ,g_p  \mid g  \mid a_0 , \dots , a_q )\\ 
&+& \sum _{1 \le i \le p} (-1)^i (g_1,\dots, g_i^{(0)}, g_i^{(1)},\dots, g_p \mid g  \mid a_0 , \dots , a_q)\\ 
&+& (-1)^{(p+1)} ( g_1, \dots ,g_p ,S^{-1}(g^{(2)})(a_0^{\bar{(1)}} \dots a_q^{\bar{(1)}})g^{(0)}\mid g^{(1)} \mid a_0^{\bar{(0)}} , \dots , a_q^{\bar{(0)}} )\\
&=& \boldsymbol{\delta} (g_1, \dots,g_p \mid g \mid a_0,\dots,a_q).
\end{eqnarray*}
So we have,
\begin{eqnarray*}
& &\bar{\mathfrak{d}}_0 (g_1, \dots ,g_p \mid g \mid a_0, \dots , a_q) = (1,g_1 , \dots , g_p \mid g \mid a_0 , \dots , a_q),\\
& &\bar{\mathfrak{d}}_i (g_1, \dots ,g_p \mid g \mid a_0, \dots , a_q) = (g_1 , \dots , g_i^{(0)}, g_i^{(1)} ,\dots , g_p \mid g \mid a_0 , \dots , a_q)\;\; 1 \le i \le p,\\
& &\bar{\mathfrak{d}}_{p+1} (g_1, \dots ,g_p \mid g \mid a_0, \dots , a_q) = (g_1 , \dots , g_{p} \mid \Delta( g \mid a_0 , \dots ,  a_q)).\\
\end{eqnarray*}
$ \hspace{0cm} \bullet \quad $ To compute $\bar{\mathfrak{t}}= \boldsymbol{\beta} \bar{\tau} \boldsymbol{\gamma}$, we see that,
\begin{eqnarray*}
& &\bar{\mathfrak{t}} (g_1 , \dots , g_p \mid g \mid a_0, \dots , a_q) =
\boldsymbol{\beta} \bar{\tau} (g^{(0)} , g^{(1)} \cdot (g_1, \dots , g_p) \mid a_0, \dots ,a_q)\\
&=& \boldsymbol{\beta} (g^{(1)} \cdot (g_1, \dots ,g_p), (a_0^{\bar{(1)}} \dots a_q^{\bar{(1)}}) g^{(0)} 
\mid a_0^{\bar{(0)}}, \dots , a_q^{\bar{(0)}})\\
&=& \boldsymbol{\beta}(g^{(1)} g_1, \dots ,g^{(p)} g_p,(a_0^{\bar{(1)}} \dots a_q^{\bar{(1)}})g^{(0)} \mid a_0^{\bar{(0)}} , \dots ,a_q^{\bar{(0)}})\\
&=& (S^{-1}(g^{(2)} g_1^{(1)}) \cdot (g^{(3)}g_2, \dots , g^{(p+1)} g_p , (a_0^{\bar{(1)}} \dots a_q^{\bar{(1)}}) g^{(0)}) 
\mid g^{(1)} g_1^{(0)} \mid a_0^{\bar{(0)}} , \dots , a_q^{\bar{(0)}})\\
&=& (S^{-1}(g_1^{(1)}) \cdot (g_2, \dots , g_p , S^{-1}(g^{(2)})(a_0^{\bar{(1)}} \dots a_q^{\bar{(1)}})g^{(0)}) \mid g^{(1)} 
g_1^{(0)} \mid a_0^{\bar{(0)}} , \dots ,a_q^{\bar{(0)}}).
\end{eqnarray*}
So we conclude that
\begin{eqnarray*}
& & \bar{\mathfrak{t}} (g_1 , \dots , g_p \mid g \mid a_0, \dots , a_q) \\
&=& (S^{-1}(g_1^{(1)}) \cdot (g_2, \dots , g_p , S^{-1}(g^{(2)})(a_0^{\bar{(1)}} \dots a_q^{\bar{(1)}})g^{(0)}) \mid g^{(1)} 
g_1^{(0)} \mid a_0^{\bar{(0)}} , \dots ,a_q^{\bar{(0)}}).
\end{eqnarray*}
For  $\bar{\mathfrak{s}}_i= \boldsymbol{\beta} \bar{\sigma}_i \boldsymbol{\gamma}$, $ 0 \le i < p $, we have
\begin{eqnarray*}
& &\bar{\mathfrak{s}}_i (g_1 , \dots , g_p \mid g \mid a_0, \dots , a_q) = 
\boldsymbol{\beta} \bar{\sigma}_i (g^{(0)} , g^{(1)} \cdot (g_1, \dots , g_p) \mid a_0, \dots ,a_q)\\
&=& \boldsymbol{\beta} (g^{(0)} , g^{(1)} \cdot (g_1, \dots ,g_i,g_{i+2},\dots, g_p) \mid a_0, \dots ,a_q) \epsilon(g_{i+1})\\
&=& (g_1 , \dots ,g_i,g_{i+2}, \dots, g_p \mid g \mid a_0, \dots , a_q) \epsilon(g_{i+1}).
\end{eqnarray*}

$\text{}\\$
$ \hspace{0.2cm} \bullet \quad $ Next we compute the operator  $\mathfrak{b} = \boldsymbol{\beta} b \boldsymbol{\gamma}$. We have 
\begin{eqnarray*}
& &\mathfrak{b}(g_1, \dots , g_p \mid g \mid a_0 , \dots , a_q) = \boldsymbol{\beta} b (g^{(0)},g^{(1)} \cdot
(g_1, \dots ,g_p) \mid a_0 , \dots , a_q )\\
&=& \boldsymbol{\beta} (\sum_{0 \le i \le q} (-1)^i (g^{(0)},g^{(1)} \cdot (g_1, \dots ,g_p) \mid a_0 , \dots ,a_i^{(0)},a_i^{(1)},\dots, a_q )\\
&+& (-1)^{q+1}( S^{-1}(a_0^{(0) \bar{(1)}}) \cdot (g^{(0)},g^{(1)} \cdot (g_1, \dots ,g_p) \mid a_0^{(1)}, a_0, \dots , a_q,a_0^{(0) \bar{(0)}})\\
&=& (\sum_{0 \le i \le q} (-1)^i ((S^{-1}(g^{(1)})g^{(2)}) \cdot (g_1, \dots ,g_p) \mid g^{(0)} \mid a_0 , \dots ,a_i^{(0)},a_i^{(1)},\dots, a_q )\\
&+& (-1)^{q+1}(( S^{-1}( g^{(1)}) a_0^{(0) \bar{(2)}} S^{-1}(a_0^{(0) \bar{(3)}})) \cdot (g^{(2)} \cdot
(g_1, \dots ,g_p)) \mid S^{-1}(a_0^{(0) \bar{(1)}}) g^{(0)} \\ 
& & \hspace{9.5cm} \mid a_0^{(1)}, a_0, \dots , a_q,a_0^{(0) \bar{(0)}})\\
\end{eqnarray*}
\begin{eqnarray*}
&=& (\sum_{0 \le i \le q} (-1)^i (g_1, \dots ,g_p \mid g  \mid a_0 , \dots ,a_i^{(0)},a_i^{(1)},\dots, a_q )\\
&+& (-1)^{q+1}
(g_1, \dots ,g_p \mid S^{-1}(a_0^{(0) \bar{(1)}}) g \mid a_0^{(1)}, a_0, \dots , a_q,a_0^{(0) \bar{(0)}}).
\end{eqnarray*}

So we conclude that,
\begin{eqnarray*}
& &\mathfrak{d}_i(g_1,\dots,g_p \mid g \mid a_0,\dots,a_q)=(g_1 , \dots , g_p \mid g \mid a_0, \dots,a_{i}^{(0)}, a_{i}^{(1)},\dots,a_{q}),\;\; 0 \le i \le q, \\
& &\mathfrak{d}_{q+1} (g_1,\dots,g_p \mid g \mid a_0,\dots,a_q)=(g_1, \dots ,g_p \mid S^{-1}(a_0^{(0) \bar{(1)}}) g \mid a_0^{(1)}, a_0, \dots , a_q,a_0^{(0) \bar{(0)}}). 
\end{eqnarray*}

$ \hspace{0.2cm} \bullet \quad $ We consider  $\mathfrak{s}_i = \boldsymbol{\beta} \sigma_i \boldsymbol{\gamma}$ and $\mathfrak{t} = \boldsymbol{\beta} \tau \boldsymbol{\gamma}$. \\

$\mathfrak{s}_i (g_1,\dots,g_p \mid g \mid a_0, \dots ,a_q) $
\begin{eqnarray*}
&=& \boldsymbol{\beta} (g^{(0)},g^{(1)} \cdot (g_1, \dots ,g_p) \mid a_0 , \dots ,a_i,a_{i+2},\dots, a_q ) \epsilon{(a_{i+1})} \\
&=& (g_1, \dots ,g_p \mid g \mid a_0 , \dots ,a_i,a_{i+2},\dots, a_q ) \epsilon{(a_{i+1})}. \\
\end{eqnarray*}

Also, for
$\mathfrak{t}= \boldsymbol{\beta} \tau \boldsymbol{\gamma} $, we have
\begin{eqnarray*}
& & \mathfrak{t} (g_1, \dots , g_p \mid g \mid a_0. \dots , a_q ) \\
&=& \boldsymbol{\beta} \tau (g^{(0)},g^{(1)} \cdot (g_1, \dots ,g_p) \mid a_0 , \dots , a_q )\\
&=& \boldsymbol{\beta} (S^{-1}(a_0^{\bar{(1)}}) \cdot (g^{(0)},g^{(1)} \cdot (g_1, \dots ,g_p)) \mid a_1 , \dots , a_q ,a_0^{\bar{(0)}}) )\\    
&=& \boldsymbol{\beta} (S^{-1}(a_0^{\bar{(1)}}) g^{(0)},S^{-1}(a_0^{\bar{(2)}}) \cdot (g^{(1)} \cdot (g_1, \dots ,g_p)) \mid a_1 , \dots , a_q ,a_0^{\bar{(0)}}) )\\    
&=& ((S^{-1}(g^{(1)}) a_0^{\bar{(2)}} S^{-1}(a_0^{\bar{(3)}}) ) \cdot (g_1, \dots ,g_p)) \mid S^{-1}(a_0^{\bar{(1)}}) g^{(0)} \mid a_1 , \dots , a_q ,a_0^{\bar{(0)}}) )\\    
&=& ((S^{-1}(g^{(1)}) \cdot (g_1, \dots ,g_p)) \mid S^{-1}(a_0^{\bar{(1)}}) g^{(0)} \mid a_1 , \dots , a_q ,a_0^{\bar{(0)}}) ).\\    
\end{eqnarray*}

By the above computations, we can state the following theorems.
\begin{theorem} \label{th:iso1}
The complex $(\mathsf{C} (C \natural^{} \mathcal{H})_{p,q} \boxtimes \mathsf{W}, b + \mathbf{u} B , \bar{b} + \mathbf{u} \bar{B}) $ 
is isomorphic to the complex 
$(\mathsf{C}^p(\mathcal{H}, \mathsf{C}^q(C^{\natural^{}}_{\mathcal{H}} \boxtimes \mathsf{W}), \mathfrak{b} + \mathbf{u} \mathfrak{B} , \bar{\mathfrak{b}} + \mathbf{u} \bar{\mathfrak{B}}), $
where $\mathfrak{\bar{b}}$ is the Hopf-comodule coboundary on $\mathsf{C}^{q}(C^{\natural^{}}_{\mathcal{H}})$. 
\end{theorem}
\begin{theorem} \label{th:iso2}
For  $p \ge 0$, $H^p(\mathcal{H},\mathsf{C}^{q} (C^{\natural^{}}_{\mathcal{H}}))$, with the operators defined as follows,
are cocyclic modules, where
$H^0(\mathcal{H},\mathsf{C}^{q}(C^{\natural^{}}_{\mathcal{H}})) \cong \mathsf{C}^{\bullet}_{\mathcal{H}^{}}(C) $, 
\begin {eqnarray*} \label{eq:equ}
& &\mathfrak{t} (g_1, \dots ,g_p \mid g \mid a_0 , \dots , a_q )=
((S^{-1}(g^{(1)}) \cdot (g_1, \dots ,g_p)) \mid S^{-1}(a_0^{\bar{(1)}}) g^{(0)} \mid a_1 , \dots , a_q ,a_0^{\bar{(0)}}) ),\\    
& &\mathfrak{d}_i(g_1,\dots,g_p \mid g \mid a_0,\dots,a_q)=(g_1 , \dots , g_p \mid g \mid a_0, \dots,a_{i}^{(0)}, a_{i}^{(1)},\dots,a_{q}), 0 \le i \le q, \\
& &\mathfrak{d}_{q+1} (g_1,\dots,g_p \mid g \mid a_0,\dots,a_q)=(g_1, \dots ,g_p \mid S^{-1}(a_0^{(0) \bar{(1)}}) g \mid a_0^{(1)}, a_0, \dots , a_q,a_0^{(0) \bar{(0)}}),\\ 
& &\mathfrak{s}_i(g_1, \dots ,g_p \mid g \mid a_0 , \dots , a_q )=(g_1, \dots , g_p \mid g \mid  a_0 , \dots ,a_i, a_{i+2},\dots, a_q) \epsilon(a_{i+1}).
\end{eqnarray*}
\end{theorem}
\begin{proof}
We prove this for $p=0,$ the general case being similar.

The operators are well defined on the coinvariant subspaces. For example, for $\mathfrak{t}$, since $ \mathfrak{t} \bar{\mathfrak{b}} = \bar{\mathfrak{b}} \mathfrak{t} $
and  $\bar{\mathfrak{b}} (g \mid a_0, \dots , a_n ) = \boldsymbol{\Delta} (g \mid a_0, \dots , a_n ) - (1 \mid g \mid a_0, \dots , a_n ),$ if $ \mathfrak{t} \bar{\mathfrak{b}} = 0$
then $\bar{\mathfrak{b}} \mathfrak{t} = 0$ .\\
To show for example $ \mathfrak{t} ^{n+1} = 1,$ we see that in the coinvariant subspace we have,
\begin{eqnarray*}
(1 \mid g \mid a_0, \dots , a_n ) \equiv ( S^{-1}(g^{(2)})(a_0^{\bar{(1)}} \dots a_q^{\bar{(1)}}) g^{(0)} \mid g^{(1)} \mid a_0^{\bar{(0)}}, \dots ,a_q^{\bar{(0)}}).\\ 
\end{eqnarray*}
So we have, 
\begin{eqnarray*}
& &\bar{\mathfrak{s}}_0 \circ \bar{\mathfrak{t}} (1 \mid g \mid a_0, \dots , a_n ) = \bar{\mathfrak{s}}_0 \circ \bar{\mathfrak{t}} ( S^{-1}(g^{(2)})(a_0^{\bar{(1)}} \dots a_q^{\bar{(1)}}) g^{(0)} \mid g^{(1)} \mid a_0^{\bar{(0)}}, \dots ,a_q^{\bar{(0)}})\\ 
&=&\bar{\mathfrak{s}}_0 ( S^{-1}(g^{(2)})(a_0^{\bar{(1)}} \dots a_q^{\bar{(1)}}) g^{(0)} \mid g^{(1)} \mid a_0^{\bar{(0)}}, \dots ,a_q^{\bar{(0)}})\\ 
&=& (  g^{(1)} \mid a_0^{\bar{(0)}}, \dots ,a_q^{\bar{(0)}}) \epsilon (S^{-1}(g^{(2)})(a_0^{\bar{(1)}} \dots a_q^{\bar{(1)}})g^{(0)})\\ 
&=& (g \mid a_0, \dots , a_n ). 
\end{eqnarray*}
Therefore,
\begin{eqnarray*}
& &\bar{\mathfrak{s}}_0 \circ \bar{\mathfrak{t}} ( S^{-1}(g^{(2)})(a_0^{\bar{(1)}} \dots a_q^{\bar{(1)}}) g^{(0)} \mid g^{(1)} \mid a_0^{\bar{(0)}}, \dots ,a_q^{\bar{(0)}})\\
&=&\bar{\mathfrak{s}}_0 (S^{-1}(g^{(3)}) (a_0^{\bar{(1)}} \dots a_n^{\bar{(1)}} ) g^{(1)} \mid g^{(2)} 
S^{-1}(g^{(4)}) ( a_0^{\bar{(2)}} \dots a_n^{\bar{(2)}})g^{(0)} \mid a_0^{\bar{(0)}},\dots,a_0^{\bar{(0)}})\\
&=&(  g^{(2)} S^{-1}(g^{(4)}) ( a_0^{\bar{(2)}} \dots a_n^{\bar{(2)}})g^{(0)} \mid a_0^{\bar{(0)}},\dots,a_0^{\bar{(0)}})
\epsilon(S^{-1}(g^{(3)}) (a_0^{\bar{(1)}} \dots a_n^{\bar{(1)}} ) g^{(1)})\\
&=&( ( a_0^{\bar{(1)}} \dots a_n^{\bar{(1)}})g  \mid a_0^{\bar{(0)}},\dots,a_0^{\bar{(0)}}).\\
\end{eqnarray*}
Therefore, in the coinvariant subspace we have \[(g \mid a_0, \dots , a_n ) 
=( ( a_0^{\bar{(1)}} \dots a_n^{\bar{(1)}})g  \mid a_0^{\bar{(0)}},\dots,a_n^{\bar{(0)}}), \]
and then we conclude that 
\begin{eqnarray*}
& &\mathfrak{t}^{n+1}(g \mid a_0, \dots , a_n )\\ 
&=&\mathfrak{t}^{n+1}( ( a_0^{\bar{(1)}} \dots a_n^{\bar{(1)}})g  \mid a_0^{\bar{(0)}},\dots,a_n^{\bar{(0)}}) \\
&=& (S^{-1}(a_0^{\bar{(1)}} \dots a_n^{\bar{(1)}})( a_0^{\bar{(2)}} \dots a_n^{\bar{(2)}})g  \mid a_0^{\bar{(0)}},\dots,a_n^{\bar{(0)}}) \\
&=& (g \mid a_0, \dots , a_n ).
\end{eqnarray*}
Also when $p=0$, then 
\begin{eqnarray*}
& &\boldsymbol{\delta_0}(g \mid a_0, \dots ,a_n) \\
&=& (1 \mid g \mid a_0, \dots , a_n) - \boldsymbol{\Delta} (g \mid a_0, \dots , a_n )\\ 
&=& (1 \mid g \mid a_0, \dots , a_n) - ( S^{-1}(g^{(2)})(a_0^{\bar{(1)}} \dots a_q^{\bar{(1)}}) g^{(0)} \mid g^{(1)} \mid a_0^{\bar{(0)}}, \dots ,a_q^{\bar{(0)}}). 
\end{eqnarray*}
So  $\ker \boldsymbol{\delta_0}$ is the coinvariant subspace.
\end{proof}

To compute the cohomologies of the mixed complex $(Tot^{\bullet}(\mathsf{C}(C \natural^{} \mathcal{H}), b + \bar{b} + \mathbf{u} (B + \bar{B})),$
we filter the complex $\mathsf{C}^{\bullet} (C \natural^{} \mathcal{H})$ by the subspaces \[
\mathsf{F}^i_{pq} Tot^{\bullet}(\mathsf{C}(C \natural^{} \mathcal{H} ) \boxtimes \mathsf{W}) = \sum_{q \le i} (\mathcal{H}^{\otimes(p+1)} \otimes C^{\otimes(q+1)}) \boxtimes \mathsf{W}. \]
If we separate the operator $b+ \bar{b} + \mathbf{u}(B + T \bar{B})$ as $\bar{b}+(b+  \mathbf{u}B) + \mathbf{u}T \bar{B}$. From Theorem(\ref{th:iso1}) and (\ref{th:iso2}) we can deduce the 
following theorem. 

\begin{theorem}
The $\mathsf{E}^0$-term of the spectral sequence is isomorphic to the complex \[ \mathsf{E}^0_{pq}=(\mathsf{C}^p(\mathcal{H},\mathsf{C}^q(C ^{\natural^{}}_{\mathcal{H}}) \boxtimes \mathsf{W} ),\boldsymbol{\delta)}, \]
and the $\mathsf{E}^1$-term is  \[
\mathsf{E}^1_{pq}= (H^p( \mathcal{H} , \mathsf{C}^q(C ^{\natural^{}}_{\mathcal{H}}) \boxtimes \mathsf{W} ) , \mathfrak{b} + \mathbf{u} \mathfrak{B})). \]
The $\mathsf{E}^2$-term of the spectral sequence is
\[ \mathsf{E}^2_{pq} = HC^q( H^p(\mathcal{H}, \mathsf{C}^q(C ^{\natural^{}}_{\mathcal{H}}));\mathsf{W}), \]
the cyclic cohomologies of the cocyclic modules $H^p(\mathcal{H},\mathsf{C}^q(C ^{\natural^{}}_{\mathcal{H}}))$ with coefficieints in $\mathsf{W}.$ 
\end{theorem}

We give an application of the above spectral sequence. Let $k$ be a field. Recall that a Hopf algebra $\mathcal{H}$ over
$k$ is called cosemisimple if $\mathcal{H}$ is cosemisimple as a coalgebra, that is, every left $\mathcal{H}$-comodule is
completely reducible~\cite{sw69}. One knows that a Hopf algebra is cosemisimple if and only if there exists a left integral
$x \in \mathcal{H^{\star}}$ with $x(1)=1$~(\cite{sw69},Theorem 14.0.3). It is easy to see that if $\mathcal{H}$ is 
cosemisimple and $M$ is an $\mathcal{H}$-bicomodule, then the coalgebra (Hochschild) cohomology groups satisfy
$H^i(\mathcal{H},M)=0$ for $i > 0,$ and $H^0(\mathcal{H},M)= M^{co\mathcal{H}},$ the subspace of coinvariants of 
the bicomodule $M$. In fact, we have the following homotopy 
operator $h: \mathcal{H}^{\otimes n} \otimes M \rightarrow \mathcal{H}^{\otimes(n-1)} \otimes M, n \ge 1,$
\[
h(g_1,\dots,g_n,m)=x(g_1)(g_2,\dots,g_n,m).
\]
One can check theat $\delta h + h \delta=id.$ Note that the antipode of $\mathcal{H}$ is bijective if $\mathcal{H}$ is 
cosemisimple.
\begin{prop}
Let $\mathcal{H}$ be a cosemisimple Hopf algebra. Then there is a natural isomorphism of cyclic and Hochschild cohomology
groups 
\[ HC^{\bullet}(C  \!> \! \blacktriangleleft \!\mathcal{H}) \simeq HC^{\bullet}(\mathsf{C}^{\bullet}_{\mathcal{H}}(C)),\]
\[ HH^{\bullet}(C  \!> \! \blacktriangleleft \!\mathcal{H}) \simeq HH^{\bullet}(\mathsf{C}^{\bullet}_{\mathcal{H}}(C)).\]
\end{prop}
\begin{proof}
Since $\mathcal{H}$ is cosemisimple, we have $\mathsf{E}_{p,q}^2=0$ for $p > 0$ and the spectral sequence collapses. The first column
of $\mathsf{E}^2$ is exactly $H^{\bullet}(\mathcal{H}, C^{\natural}_{\mathcal{H}})= \mathsf{C}^{\bullet}_{\mathcal{H}}(C)$.
\end{proof}


\end{document}